\documentclass[11pt]{amsart}
\usepackage{graphicx, color}
\usepackage{amscd}
\usepackage{amsmath}
\usepackage{amsfonts}
\usepackage{amssymb}
\usepackage{mathrsfs}

\newtheorem{theorem}{Theorem}

\newtheorem{remark}[theorem]{Remark}

\numberwithin{equation}{section}

\renewcommand{\leq}{\leqslant}

\renewcommand{\geq}{\geqslant}
\begin{document}
\title[]{\textbf{%
Symmetric solutions for a partial differential elliptic equation that arises
in stochastic production planning with production constraints}}
\author[D.-P. Covei]{Dragos-Patru Covei\\Department of Applied Mathematics, The Bucharest University of
Economic Studies\\
Piata Romana, 1st district, postal code: 010374, postal office: 22, Romania}
\address{}
\email{\texttt{coveidragos@yahoo.com}}

\begin{abstract}
In this article we consider the question of the existence of positive
symmetric solutions to the problems of the following type $\Delta u=a\left(
\left\vert x\right\vert \right) h\left( u\right) +b\left( \left\vert
x\right\vert \right) g\left( u\right) $ for $x\in \mathbb{R}^{N}$, which are
called entire large solutions. Here $N\geq 3$, and we assume that $a$ and $b$
are nonnegative continuous spherically symmetric functions on $\mathbb{R}%
^{N} $. We extend results previously obtained for special cases of $h$ and $%
g $ and we will describe a real-world model in which such problems may arise.
\end{abstract}

\maketitle




\section{Introduction}

The problem of establishing the existence of spherically symmetric solutions
(i.e. $u\left( x\right) =u\left( r\right) $ where $r:=\left\vert
x\right\vert $ is the Euclidean norm) for the partial differential equation 
\begin{equation}
\Delta u=a\left( \left\vert x\right\vert \right) h\left( u\right) +b\left(
\left\vert x\right\vert \right) g\left( u\right) \text{ for }x\in \mathbb{R}%
^{N}\text{ (}N\geq 3\text{)},  \label{11}
\end{equation}%
is extensively studied both from the theoretical point of view and from
modeling of various phenomena in the real world. From a theoretical point of
view, using conditions on the potential functions 
\begin{equation*}
a,b\in C_{loc}^{0,\beta }\left( \mathbb{R}^{N},\left[ 0,\infty \right)
\right) \text{ for some }\beta \in \left( 0,1\right) \text{,}
\end{equation*}%
and on the nonlinearities 
\begin{equation*}
h,g\in C\left( \left[ 0,\infty \right) ,\left[ 0,\infty \right) \right) \cap
C^{1}\left( \left( 0,\infty \right) ,\left[ 0,\infty \right) \right) ,
\end{equation*}%
the study of the existence solutions for the problem (\ref{11}) was well
argued, for the case where $h$ and $g$ are nondecreasing with $h\left(
0\right) =g\left( 0\right) =0$ and $h\left( s\right) g\left( s\right) >0$
for all $s>0$, in the article of Lair and Mohammed [\cite{LM}]. Inspired by
the works of Alvarez [\cite{AL}], Bensoussan, Sethi, Vickson and Derzko [%
\cite{B}], Du and Guo [\cite{DG}], Lair and Mohammed [\cite{LM}], Lasry and
Lions [\cite{LS}]\ and Porretta [\cite{P}] here we assume that the
nonlinearities $h$ and $g$ belongs to a new class of functions:

h1) \quad $h:\left[ 0,\infty \right) \rightarrow \left[ 0,\infty \right) $
is continuous, nondecreasing with $h\left( 0\right) =0$ and $h\left(
s\right) >0$ for all $s\in \left( 0,\infty \right) $;

g1)\quad $g:\left( 0,\infty \right) \rightarrow \mathbb{R}$ is continuous
and there exists $s_{0}\in \left( 0,\infty \right) $ such that $g$ is
non-decreasing for all $s\in \left( s_{0},\infty \right) $, $%
\lim_{s\rightarrow 0+}g\left( s\right) =g\left( s_{0}\right) =0$, $g\left(
s\right) <0$ for all $s\in \left( 0,s_{0}\right) $ and $g\left( s\right) >0$
for all $s\in \left( s_{0},\infty \right) $;

O)\quad for $s_{0}\in \left( 0,\infty \right) $, which exists from the
condition g1), we have%
\begin{equation*}
H\left( \infty \right) =\infty \text{ for all }u_{0}>s_{0},
\end{equation*}%
where%
\begin{equation*}
H\left( \infty \right) :=\lim_{t\rightarrow \infty }H\left( t\right) \text{, 
}H\left( t\right) :=\int_{u_{0}}^{t}\frac{1}{h\left( t\right) +g\left(
t\right) }dt\text{, }t\geq u_{0}\text{,}
\end{equation*}%
and we are interested in the following goals for the problem (\ref{11}):
establishing a result of existence of spherically symmetric solutions for (%
\ref{11}), determining the asymptotic behavior of solutions for (\ref{11})
and last but not least is to write a description of a model from real-world
where such problems (\ref{11}) might arise.

More exactly, that we aim to emphasize on problems (\ref{11}) from the
theoretical point of view, is synthesized in:

\begin{theorem}
\label{th1} We assume h1), g1) and O) hold. If $a$ and $b$ are nonnegative
continuous spherically symmetric functions on $\mathbb{R}^{N}$, then 
\begin{equation*}
\left\{ 
\begin{array}{l}
u\left( 0\right) =u_{0}\geq s_{0}, \\ 
u_{n}\left( r\right) =u_{0}+\int_{0}^{r}t^{1-N}\int_{0}^{t}s^{N-1}(a\left(
s\right) h\left( u_{n-1}\left( s\right) \right) +b\left( s\right) g\left(
u_{n-1}\left( s\right) \right) )dsdt\text{, }r\geq 0%
\end{array}%
\right.
\end{equation*}%
converges locally uniformly to a positive sperically symmetric function $%
u\in C^{2}\left[ 0,\infty \right) $ and $u$ is a solution of the equation (%
\ref{11}) such that $u\left( 0\right) =u_{0}$.

In addition, for any $c_{0}>0$ we have the estimates:

A1)\quad there exists $\overline{C}_{1}>u_{0}$ such that $u_{0}\leq u\left(
r\right) \leq \overline{C}_{1}$, for all $r\in \left[ 0,c_{0}\right] ;$

A2)\quad there exists $\overline{C}_{2}>0$ such that $0\leq u^{\prime
}\left( r\right) \leq \overline{C}_{2}\left( r+1\right) $, for all $r\in %
\left[ 0,c_{0}\right] .$

Moreover, with the following notations 
\begin{eqnarray*}
m\left( s\right) &:&=\left\{ 
\begin{array}{lll}
a\left( s\right) h\left( u_{0}\right) & \text{if} & g\left( u_{0}\right) =0%
\text{,} \\ 
\left( \left( a\left( s\right) +b\left( s\right) \right) \min \left\{
h\left( u_{0}\right) ,g\left( u_{0}\right) \right\} \right) & \text{if} & 
g\left( u_{0}\right) \neq 0\text{,}%
\end{array}%
\right. \\
\underline{P}\left( r\right)
&:&=\int_{0}^{r}t^{1-N}\int_{0}^{t}s^{N-1}m\left( s\right) dsdt\text{, }%
\underline{P}\left( \infty \right) :=\lim_{r\rightarrow \infty }\underline{P}%
\left( r\right) \text{,} \\
\overline{P}\left( r\right) &:&=\int_{0}^{r}t^{1-N}\int_{0}^{t}s^{N-1}\left(
a\left( s\right) +b\left( s\right) \right) dsdt\text{, }\overline{P}\left(
\infty \right) :=\lim_{r\rightarrow \infty }\overline{P}\left( r\right) 
\text{,}
\end{eqnarray*}%
if $\underline{P}\left( \infty \right) =\infty $ then $\underset{%
r\rightarrow \infty }{\lim }u\left( r\right) =\infty $ (i.e. $u$ is large),
and $\underset{r\rightarrow \infty }{\lim }u\left( r\right) <\infty $
provided $\overline{P}\left( \infty \right) <\infty $ (i.e. $u$ is bounded).
\end{theorem}

\begin{remark}
\label{rem}If $a$ and $b$\textit{\ are nondecreasing functions }then the
solution $u$ obtained in the Theorem \ref{th1} is convex.
\end{remark}

In the following we will structure the paper in two sections. In the Section
2 we will prove the Theorem \ref{th1} and Remark \ref{rem}. For this, we use
a different approach to that of Lasry and Lions [\cite{LS}] or Alvarez [\cite%
{AL}] at least from the following aspects: approximating the solution for
the considered problem (\ref{11}), the direct study of the solved problem
without call to auxiliary results, imposing another conditions on the
boundary and last but not least by establishing the monotony of the
solution. The difficulty that arises from these aspects is given by the
influence of the sign of the function $g$. In the Section 3 we will describe
a real-world model in which such problems may arise. In addition, new
directions of study are open to research.

\section{The Proof of Theorem}

The main objective, is to establish the existence of solutions for the
second order differential equation 
\begin{equation}
\left\{ 
\begin{array}{l}
\left( r^{N-1}u^{\prime }\left( r\right) \right) ^{\prime }=r^{N-1}(a\left(
r\right) h\left( u\left( r\right) \right) +b\left( r\right) g\left( u\left(
r\right) \right) )\text{, }r>0, \\ 
u\left( 0\right) =u_{0}\in \left[ s_{0},\infty \right) \text{, }u^{\prime
}\left( 0\right) =0\text{.}%
\end{array}%
\right.  \label{ss1}
\end{equation}%
A common argument is to rewrite (\ref{ss1}) in the integral form 
\begin{equation}
u\left( r\right) =u_{0}+\int_{0}^{r}t^{1-N}\int_{0}^{t}s^{N-1}(a\left(
s\right) h\left( u\left( s\right) \right) +b\left( s\right) g\left( u\left(
s\right) \right) )dsdt.  \label{ss}
\end{equation}%
So, setting a solution for (\ref{ss}) we get the solution for the proposed
problem (\ref{11}). We define the sequence of functions $\left\{ u_{n}\left(
r\right) \right\} _{n\in \mathbb{N}}$ in $\left[ 0,\infty \right) $
iteratively, as follows: 
\begin{equation}
\left\{ 
\begin{array}{l}
u_{0}\left( r\right) :=u_{0}\geq s_{0}, \\ 
u_{n}\left( r\right) =a+\int_{0}^{r}t^{1-N}\int_{0}^{t}s^{N-1}(a\left(
s\right) h\left( u_{n-1}\left( s\right) \right) +b\left( s\right) g\left(
u_{n-1}\left( s\right) \right) dsdt\text{, }r\geq 0.%
\end{array}%
\right.  \label{recs}
\end{equation}%
We note that for any $r\in \left[ 0,\infty \right) $ the sequence of
functions $\left\{ u_{n}\left( r\right) \right\} _{n\in \mathbb{N}}$ it is
monotonically increasing with respect to variable $r$ independent of the
value of $n\in \mathbb{N}$, which is a useful information in the following.
To prove the existence of the limit%
\begin{equation*}
u\left( r\right) :=\lim_{n\rightarrow \infty }u_{n}\left( r\right) ,
\end{equation*}%
we will prove that$\left\{ u_{n}\left( r\right) \right\} _{n\in \mathbb{N}}$
is a nondecreasing sequence on $\left[ 0,\infty \right) $. We use the
mathematical induction method. We note that the verification step takes
place 
\begin{eqnarray*}
u_{1}\left( r\right) &=&u_{0}+\int_{0}^{r}t^{1-N}\int_{0}^{t}s^{N-1}\left(
a\left( s\right) h\left( u_{0}\left( s\right) \right) +b\left( s\right)
g\left( u_{0}\left( s\right) \right) \right) dsdt \\
&=&u_{0}+\int_{0}^{r}t^{1-N}\int_{0}^{t}s^{N-1}\left( a\left( s\right)
h\left( u_{0}\right) +b\left( s\right) g\left( u_{0}\right) \right) dsdt \\
&\leq &u_{0}+\int_{0}^{r}t^{1-N}\int_{0}^{t}s^{N-1}\left( a\left( s\right)
h\left( u_{1}\left( s\right) \right) +b\left( s\right) g\left( u_{1}\left(
s\right) \right) \right) dsdt \\
&=&u_{2}\left( r\right) ,
\end{eqnarray*}%
and so $u_{1}\left( r\right) \leq u_{2}\left( r\right) $ for any $r\in \left[
0,\infty \right) $. We assume $u_{n-1}\left( r\right) \leq u_{n}\left(
r\right) $ for any $r\in \left[ 0,\infty \right) $ and we prove that%
\begin{equation*}
u_{n}\left( r\right) \leq u_{n+1}\left( r\right) \text{ for any }n\in 
\mathbb{N}\text{ and }r\in \left[ 0,\infty \right) .
\end{equation*}%
Indeed,%
\begin{eqnarray*}
u_{n}\left( r\right) &=&u_{0}+\int_{0}^{r}t^{1-N}\int_{0}^{t}s^{N-1}\left(
a\left( s\right) h\left( u_{n-1}\left( s\right) \right) +b\left( s\right)
g\left( u_{n-1}\left( s\right) \right) \right) dsdt \\
&\leq &u_{0}+\int_{0}^{r}t^{1-N}\int_{0}^{t}s^{N-1}\left( a\left( s\right)
h\left( u_{n}\left( s\right) \right) +b\left( s\right) g\left( u_{n}\left(
s\right) \right) \right) dsdt \\
&=&u_{n+1}\left( r\right) ,
\end{eqnarray*}%
which ends the proof.

In the following we will show that the sequence of functions $\left\{
u_{n}\left( r\right) \right\} _{n\in \mathbb{N}}$ is bounded in any compact
interval $\left[ 0,c_{0}\right] $ with $c_{0}>0$ arbitrary and independent
of $n$, property that couples with the monotony of the sequence ensures the
uniform convergence%
\begin{equation*}
u_{n}\left( r\right) \underset{r\in \left[ 0,c_{0}\right] }{\overset{\text{%
uniformly}}{\rightarrow }}u\left( r\right) \text{ as }n\rightarrow \infty .
\end{equation*}%
To prove this, we notice that%
\begin{eqnarray}
\left[ r^{N-1}\left( u_{n}\left( r\right) \right) ^{\prime }\right] ^{\prime
} &=&r^{N-1}\left( a\left( s\right) h\left( u_{n-1}\left( r\right) \right)
+b\left( s\right) g\left( u_{n-1}\left( r\right) \right) \right)  \notag \\
&\leq &r^{N-1}\left( a\left( s\right) h\left( u_{n}\left( r\right) \right)
+b\left( s\right) g\left( u_{n}\left( r\right) \right) \right)  \label{gen2}
\\
&\leq &r^{N-1}\left( a\left( r\right) +b\left( r\right) \right) \left(
h\left( u_{n}\left( r\right) \right) +g\left( u_{n}\left( r\right) \right)
\right) ,  \notag
\end{eqnarray}%
where we have used the monotony of the sequence $\left\{ u_{n}\left(
r\right) \right\} _{n\in \mathbb{N}}$. Integrating between $0$ and $r$ the
inequality (\ref{gen2}), we find 
\begin{eqnarray*}
\left( u_{n}\left( r\right) \right) ^{\prime } &\leq
&r^{1-N}\int_{0}^{r}t^{N-1}\left( a\left( t\right) +b\left( t\right) \right)
\left( h\left( u_{n}\left( t\right) \right) +g\left( u_{n}\left( t\right)
\right) \right) dt \\
&\leq &r^{1-N}\left( h\left( u_{n}\left( r\right) \right) +g\left(
u_{n}\left( r\right) \right) \right) \int_{0}^{r}t^{N-1}\left( a\left(
t\right) +b\left( t\right) \right) dt,
\end{eqnarray*}%
and after rearrangement 
\begin{equation}
\frac{\left( u_{n}\left( r\right) \right) ^{\prime }}{h\left( u_{n}\left(
r\right) \right) +g\left( u_{n}\left( r\right) \right) }\leq
r^{1-N}\int_{0}^{r}t^{N-1}\left( a\left( t\right) +b\left( t\right) \right)
dt.  \label{mat2}
\end{equation}%
A new integration from $0$ to $r$, in the inequality (\ref{mat2}), leads to%
\begin{equation*}
\int_{u_{0}}^{u_{n}\left( r\right) }\frac{1}{h\left( t\right) +g\left(
t\right) }dt\leq \overline{P}\left( r\right) ,\text{ }\overline{P}^{\prime
}\left( r\right) :=r^{1-N}\int_{0}^{r}t^{N-1}\left( a\left( t\right)
+b\left( t\right) \right) dt,
\end{equation*}%
expression that can be written such%
\begin{equation}
H\left( u_{n}\left( r\right) \right) \leq \overline{P}\left( r\right) \text{%
, }H\left( s\right) :=\int_{u_{0}}^{s}\frac{1}{h\left( t\right) +g\left(
t\right) }dt\text{, }s>u_{0}.  \label{ints}
\end{equation}%
Obviously the function $H:\left[ u_{0},\infty \right) \rightarrow \left[
0,\infty \right) $ is bijective and strictly increasing for any $s\in \left[
u_{0},\infty \right) $, properties that which are transmited to the inverse
function $H^{-1}$. Applying the inverse $H^{-1}$ to inequality (\ref{ints})
we obtain%
\begin{equation}
u_{n}\left( r\right) \leq H^{-1}\left( \overline{P}\left( r\right) \right) .
\label{int}
\end{equation}%
Recapitulate, it is noticed that%
\begin{equation}
u_{n}\left( r\right) \leq u_{n+1}\left( r\right) \text{ for any }n\in 
\mathbb{N},\text{ }r\in \left[ 0,\infty \right) ,  \label{d1}
\end{equation}%
and 
\begin{equation}
u_{n}\left( r\right) \leq u_{n}\left( c_{0}\right) \leq \overline{C}%
_{1}:=H^{-1}\left( \overline{P}\left( c_{0}\right) \right) <\infty \text{ }%
\forall n\in \mathbb{N},\text{ }\forall r\in \left[ 0,c_{0}\right] ,
\label{d2}
\end{equation}%
properties justifying the applicability of Dini's theorem and implicitly
establishing that 
\begin{equation*}
u_{n}\left( r\right) \underset{r\in \left[ 0,c_{0}\right] }{\overset{\text{%
uniformly}}{\rightarrow }}u\left( r\right) \text{ for }n\rightarrow \infty 
\text{ and }c_{0}>0\text{ arbitrary.}
\end{equation*}%
So, letting to the limit when $n\rightarrow \infty $ in (\ref{recs}) we
obtain the existence of a $c_{0}:=c_{0}(u_{0})>0$ (maximal extreme to the
right) for the existence maximal interval $\left( 0,c_{0}\right) $ of
solutions for (\ref{ss1}) and an $u\left( r\right) :=u_{u_{0}}\left(
r\right) \in C^{2}\left( 0,c_{0}\right) \cap C^{1}\left( \left[
0,c_{0}\right) \right) $ solution for the problem (\ref{ss1}) in $\left(
0,c_{0}\right) $. We prove that $u\left( r\right) $ exists in $\left(
0,\infty \right) $ which is reduced to showing that $c_{0}=\infty $. Assume
the contrary, that $c_{0}\in \left( 0,\infty \right) $. A simple argument
proves that $\underset{\underset{r<c_{0}}{r\rightarrow c_{0}}}{\lim }u\left(
r\right) =\infty $. Appealing to the inequality 
\begin{equation*}
\int_{u_{0}}^{u\left( r\right) }\frac{1}{h\left( t\right) +g\left( t\right) }%
dt\leq \overline{P}\left( r\right) ,
\end{equation*}%
by letting to the limit, we get 
\begin{eqnarray*}
\infty &=&\underset{\underset{r<c_{0}}{r\rightarrow c_{0}}}{\lim }%
\int_{u_{0}}^{u\left( r\right) }\frac{1}{h\left( t\right) +g\left( t\right) }%
dt \\
&=&\underset{\underset{r<c_{0}}{r\rightarrow c_{0}}}{\lim }%
\int_{u_{0}}^{\infty }\frac{1}{h\left( t\right) +g\left( t\right) }dt\leq 
\underset{\underset{r<c_{0}}{r\rightarrow c_{0}}}{\lim }\overline{P}\left(
r\right) <\infty ,
\end{eqnarray*}%
and then a contradiction. We have proved that 
\begin{equation*}
u_{n}\left( r\right) \underset{r\in \left[ 0,\infty \right) }{\overset{\text{%
locally uniformly}}{\rightarrow }}u\left( r\right) \text{ for }n\rightarrow
\infty \text{,}
\end{equation*}%
and as a consequence $u\left( r\right) $ satisfy%
\begin{equation*}
\left\{ 
\begin{array}{l}
u\left( 0\right) =u_{0} \\ 
u\left( r\right) =u_{0}+\int_{0}^{r}y^{1-N}\int_{0}^{y}t^{N-1}\left( a\left(
t\right) h\left( u\left( t\right) \right) +b\left( t\right) g\left( u\left(
t\right) \right) \right) dtdy,\text{ }r\geq 0.%
\end{array}%
\right.
\end{equation*}%
The regularity of the solution is a classic process that can be consulted in
the paper of [\cite{COV}].

We prove the announced estimates. Let $c_{0}>0$ arbitrary. The claim A1) it
is evident. Then it remains to test the existence of a parameter $\overline{C%
}_{2}>0$ such that $u^{\prime }\left( r\right) \leq \overline{C}_{2}\left(
r+1\right) $ for any $r\in \left[ 0,c_{0}\right] $. Indeed, for any $r\geq 0$%
,%
\begin{eqnarray}
u^{\prime }\left( r\right) &=&r^{1-N}\int_{0}^{r}t^{N-1}(a\left( t\right)
h\left( u\left( t\right) \right) +b\left( t\right) g\left( u\left( t\right)
\right) )dt  \notag \\
&\leq &\left( h\left( u\left( r\right) \right) +g\left( u\left( r\right)
\right) \right) \int_{0}^{r}\left( a\left( t\right) +b\left( t\right)
\right) dt  \notag \\
&\leq &\left( h\left( u\left( c_{0}\right) \right) +g\left( u\left(
c_{0}\right) \right) \right) \int_{0}^{r}\left( a\left( t\right) +b\left(
t\right) \right) dt  \label{estimat} \\
&\leq &\left\Vert a+b\right\Vert _{\infty }\left( h\left( u\left(
c_{0}\right) \right) +g\left( u\left( c_{0}\right) \right) \right)
\int_{0}^{r}dt  \notag \\
&=&\left\Vert a+b\right\Vert _{\infty }\left( h\left( u\left( c_{0}\right)
\right) +g\left( u\left( c_{0}\right) \right) \right) r\text{,}  \notag
\end{eqnarray}%
the inequalities that take place for any $r\in \left[ 0,c_{0}\right] $. As a
consequence%
\begin{equation*}
\overline{C}_{2}:=\left\Vert a+b\right\Vert _{\infty }\left( h\left( u\left(
c_{0}\right) \right) +g\left( u\left( c_{0}\right) \right) \right) ,
\end{equation*}%
check the affirmations.

Next, we chek the convexity of the solution\textbf{\ }$u$\textbf{. }Indeed,
it is clear that%
\begin{equation}
\left( r^{N-1}u^{^{\prime }}\left( r\right) \right) ^{\prime
}=r^{N-1}(a\left( r\right) h\left( u\left( r\right) \right) +b\left(
r\right) g\left( u\left( r\right) \right) ).  \label{jo1}
\end{equation}%
Integrating the equation (\ref{jo1}) from $0$ to $r>0$ we obtain%
\begin{eqnarray*}
r^{N-1}u^{^{\prime }}\left( r\right) &=&\int_{0}^{r}s^{N-1}(a\left( s\right)
h\left( u\left( s\right) \right) +b\left( s\right) g\left( u\left( s\right)
\right) )ds \\
&\leq &a\left( r\right) h\left( u\left( r\right) \right)
\int_{0}^{r}s^{N-1}ds+b\left( r\right) g\left( u\left( r\right) \right)
\int_{0}^{r}s^{N-1}ds \\
&=&a\left( r\right) h\left( u\left( r\right) \right) \frac{r^{N}}{N}+\frac{%
r^{N}}{N}b\left( r\right) g\left( u\left( r\right) \right) \\
&=&\frac{r^{N}}{N}\left( a\left( r\right) h\left( u\left( r\right) \right)
+b\left( r\right) g\left( u\left( r\right) \right) \right) ,
\end{eqnarray*}%
and, as a consequence%
\begin{equation}
\frac{u^{^{\prime }}\left( r\right) }{r}\leq \frac{1}{N}\left( a\left(
r\right) h\left( u\left( r\right) \right) +b\left( r\right) g\left( u\left(
r\right) \right) \right) \text{, }\forall r>0\text{.}  \label{jo}
\end{equation}%
On the other hand (\ref{jo1}) can be written in form 
\begin{equation}
u^{^{\prime \prime }}\left( r\right) +\left( N-1\right) \frac{u^{^{\prime
}}\left( r\right) }{r}=a\left( r\right) h\left( u\left( r\right) \right)
+b\left( r\right) g\left( u\left( r\right) \right) .  \label{fc}
\end{equation}%
substituting (\ref{jo}) into (\ref{fc}) we obtain \ 
\begin{equation}
a\left( r\right) h\left( u\left( r\right) \right) +b\left( r\right) g\left(
u\left( r\right) \right) \leq u^{^{\prime \prime }}\left( r\right) +\frac{N-1%
}{N}\left( a\left( r\right) h\left( u\left( r\right) \right) +b\left(
r\right) g\left( u\left( r\right) \right) \right) .  \label{gra}
\end{equation}%
Rearranging the inequality (\ref{gra}), we get 
\begin{equation*}
0<\left( a\left( r\right) h\left( u\left( r\right) \right) +b\left( r\right)
g\left( u\left( r\right) \right) \right) \left( 1-\frac{N-1}{N}\right) \leq
u^{^{\prime \prime }}\left( r\right) \text{, }\forall r>0,
\end{equation*}%
relation that coupled with%
\begin{equation*}
u^{^{\prime \prime }}\left( 0\right) =\frac{a\left( 0\right) h\left( u\left(
0\right) \right) +b\left( 0\right) g\left( u\left( 0\right) \right) }{N}\geq
0,
\end{equation*}%
leads to $u^{^{\prime \prime }}\left( r\right) \geq 0$ for any $r\geq 0$.

Finally, we prove the limit of the solution on the boundary. In the case $%
\underline{P}\left( \infty \right) =\infty $, we observe that 
\begin{eqnarray*}
u\left( r\right) &=&u_{0}+\int_{0}^{r}t^{1-N}\int_{0}^{t}s^{N-1}\left(
a\left( s\right) h\left( u\left( s\right) \right) +b\left( s\right) g\left(
u\left( s\right) \right) \right) dsdt \\
&\geq &u_{0}+\int_{0}^{r}t^{1-N}\int_{0}^{t}s^{N-1}\left( a\left( s\right)
h\left( u_{0}\right) +b\left( s\right) g\left( u_{0}\left( s\right) \right)
\right) dsdt \\
&\geq &u_{0}+\int_{0}^{r}t^{1-N}\int_{0}^{t}s^{N-1}m\left( s\right) dsdt,
\end{eqnarray*}%
or, considering the inequality in which are interested 
\begin{equation}
u\left( r\right) \geq u_{0}+\underline{P}\left( r\right) .  \label{bc}
\end{equation}%
Consequently, passing to the limit in (\ref{bc}) we obtain $%
\lim_{r\rightarrow \infty }u\left( r\right) =\infty $. In the case $%
\overline{P}\left( \infty \right) <\infty $, we see that 
\begin{equation*}
u\left( r\right) \leq H^{-1}\left( \overline{P}\left( r\right) \right) \leq
H^{-1}\left( \overline{P}\left( \infty \right) \right) <\infty .
\end{equation*}%
We also point that 
\begin{equation*}
\lim_{r\rightarrow \infty }\int_{0}^{r}y^{1-N}\int_{0}^{y}t^{N-1}\left(
a\left( t\right) +b\left( t\right) \right) dtdy=\frac{1}{N-2}%
\int_{0}^{\infty }r\left( a\left( r\right) +b\left( r\right) \right) dr.
\end{equation*}

\section{The model \label{sec2}}

Consider a factory producing $N$ homogeneous goods and having an inventory
warehouse. Let $p\left( t\right) =\left( p_{1}(t),...,p_{N}(t)\right) \geq 0$
represents the production rate at time $t$ (control variable) and $y\left(
t\right) =\left( y_{1}(t),...,y_{N}(t)\right) $ denote the inventory level
for production rate at time $t$ (state variable). We point that a negative
value of $y_{i}(t)$ indicates a backlogged demand for part $i$ (for example,
due to obsolescence or perishability), while a positive value is the size of
the inventory stored in the buffers. We consider $w=\left(
w_{1},...,w_{N}\right) $ a $N$-dimensional Brownian motion on a complete
probability space $(\Omega ,\mathcal{F},P)$ endowed with the natural
completed filtration $\{\mathcal{F}_{t}\}_{0\leq t\leq T}$, where $T$ is the
length of planning period, generated by the standard Wiener process $w$. We
note that the process $dw_{i}(t)$, $i=1,2,...,N$, can be formally expressed
as $z_{i}(t)dt$, where $z_{i}(t)$ is considered to be the white noise
process [\cite{A}]. For a fall in the Theorem \ref{th1} hypotheses, we
assume that all the constant demand rate at time $t$ are equal with $0$ and
we consider the inventory production control problem 
\begin{equation}
J\left( p_{1},...,p_{N}\right) :=\text{ }E\int_{0}^{\infty
}(f_{1}(p(t))+f_{2}(y(t)))e^{-\alpha t}dt,  \label{opti}
\end{equation}%
where $f_{1}(x)=f_{2}\left( x\right) =\left\vert x\right\vert ^{2}$ is the
quadratic loss function. The stochastic differential equations governing $%
y_{i}$ are 
\begin{equation}
dy_{i}\left( t\right) =p_{i}dt+\sigma _{i}dw_{i}\text{, }y_{i}\left(
0\right) =y_{i}^{0}\text{, }i=1,...,N,  \label{cons}
\end{equation}%
where $\sigma =\left( \sigma _{1},...,\sigma _{N}\right) $ is the non-zero
constant diffusion coefficient, $\alpha >0$ is the constant discount rate
and $y_{i}^{0}$ is the initial inventory level. The aim is to minimize the
stochastic production planning problem%
\begin{equation}
\inf \{J\left( p_{1},...,p_{N}\right) \left\vert p_{i}\geq 0\text{ }\forall
i=1,2,...,N\right. \}\text{, }  \label{min}
\end{equation}%
subject to the It\^{o} equation (\ref{cons}). Let $z\left( x\right) =z\left(
x_{1},...,x_{N}\right) $ denote the expected current-valued value of the
control problem (\ref{cons})-(\ref{min}) with initial value $x$ in (\ref%
{cons}) so that $z\left( y_{i}^{0}\right) $ represents to the state-equation
(\ref{cons}). In order to achieve this we apply the martingale principle: we
search for a function $U\left( x\right) $ such that the stochastic process $%
M^{c}(t)$ defined below%
\begin{equation*}
M^{c}\left( t\right) =e^{-\alpha t}U\left( y\left( t\right) \right)
-E\int_{0}^{\infty }(f_{1}(p(t))+f_{2}(y(t)))e^{-\alpha t}dt
\end{equation*}%
is supermartingale for all $p_{1}(t)\geq 0$, ... , $p_{N}(t)\geq 0$ and
martingale for the optimal control $p^{\ast }\left( t\right) =\left(
p_{1}^{\ast }(t),...,p_{N}^{\ast }(t)\right) $. If, this is achieved and the
following transversality condition holds true%
\begin{equation}
\lim_{t\rightarrow \infty }E[e^{-\alpha t}U\left( y\left( t\right) \right)
]=0\text{, }  \label{tr}
\end{equation}%
then, it can be shown that $-U\left( x\right) =z\left( x\right) $ is $C^{2}%
\left[ 0,\infty \right) $ and satisfies the associated dynamic programming
partial differential equation or Hamilton-Jacobi-Bellman equation formally
associated to the problem (\ref{cons})-(\ref{min}) 
\begin{equation}
\alpha z-\frac{\left\vert \sigma \right\vert ^{2}}{2}\Delta z-\left\vert
x\right\vert ^{2}=\inf \{p\nabla z+\left\vert p\right\vert ^{2}\left\vert
p_{i}\geq 0\text{ }\forall i=1,2,...,N\right. \},  \label{solv}
\end{equation}%
where $z:=z\left( x_{1},...,x_{N}\right) $ is the corresponding value
function. We point that the solution of this HJB equation is used to test
controller for optimality or perhaps to construct a feedback controller. In
the next, we give some ideas to solve the problem (\ref{solv}). Firstly, if $%
\frac{\partial z}{\partial x_{i}}\left( x_{1},...,x_{N}\right) \leq 0$ for
all $i=1,...,N$, since this is the case where we are interested, we observe
that%
\begin{equation*}
p\nabla z+\left\vert p\right\vert ^{2}=-\frac{1}{4}\left\vert \nabla
z\right\vert ^{2}.
\end{equation*}%
Indeed, setting%
\begin{eqnarray}
F\left( p_{1},...,p_{N}\right) &=&p\nabla z+\left\vert p\right\vert ^{2}
\label{minmin} \\
&=&p_{1}\frac{\partial z}{\partial x_{1}}\left( x_{1},...,x_{N}\right)
+...+p_{N}\frac{\partial z}{\partial x_{N}}\left( x_{1},...,x_{N}\right)
+\sum_{i=1}^{N}p_{i}^{2}  \notag
\end{eqnarray}%
we have%
\begin{equation*}
F_{p_{i}}\left( p_{1},...,p_{N}\right) =0\Longleftrightarrow \frac{\partial z%
}{\partial x_{i}}\left( x_{1},...,x_{N}\right) +2p_{i}=0
\end{equation*}%
for all $i=1,2,...,N$. Then, the critical point of the function $F$ is%
\begin{equation*}
p_{i}^{\ast }=-\frac{1}{2}\frac{\partial z}{\partial x_{i}}\left(
x_{1},...,x_{N}\right) \text{ for }i=1,...,n\text{.}
\end{equation*}%
On the other hand the Hessian matrix is%
\begin{equation*}
H_{F}\left( p_{1},...,p_{N}\right) =\left( 
\begin{array}{ccc}
2 & ... & 0 \\ 
... & ... & ... \\ 
0 & ... & 2%
\end{array}%
\right)
\end{equation*}%
which is positive definite and so 
\begin{equation*}
\left( p_{1}^{\ast },..,p_{N}^{\ast }\right) =\left( -\frac{1}{2}\frac{%
\partial z}{\partial x_{1}}\left( x_{1},...,x_{N}\right) ,...,-\frac{1}{2}%
\frac{\partial z}{\partial x_{N}}\left( x_{1},...,x_{N}\right) \right)
\end{equation*}%
is a global minimum point for (\ref{minmin}). Then, we have%
\begin{eqnarray*}
F\left( p_{1}^{\ast },...,p_{N}^{\ast }\right) &=&\left( p_{1}^{\ast
},...,p_{N}^{\ast }\right) \nabla z+\sum_{i=1}^{N}\left( p_{i}^{\ast
}\right) ^{2} \\
&=&-\frac{1}{2}\sum_{i=1}^{N}\left( \frac{\partial z}{\partial x_{i}}\left(
x_{1},...,x_{N}\right) \right) ^{2}+\frac{1}{4}\left\vert \nabla
z\right\vert ^{2} \\
&=&-\frac{1}{2}\left\vert \nabla z\right\vert ^{2}+\frac{1}{4}\left\vert
\nabla z\right\vert ^{2}=-\frac{1}{4}\left\vert \nabla z\right\vert ^{2},
\end{eqnarray*}%
so that equation (\ref{solv}) can be written as%
\begin{equation*}
\alpha z-\frac{\left\vert \sigma \right\vert ^{2}}{2}\Delta z-\left\vert
x\right\vert ^{2}=-\frac{1}{4}\left\vert \nabla z\right\vert ^{2}\text{ for }%
x\in \mathbb{R}^{N}\text{,}
\end{equation*}%
or, equivalently%
\begin{equation}
-2\left\vert \sigma \right\vert ^{2}\Delta z+\left\vert \nabla z\right\vert
^{2}+4\alpha z=4\left\vert x\right\vert ^{2}\text{ for }x\in \mathbb{R}^{N}.
\label{eq}
\end{equation}%
and, finally with the change of variable 
\begin{equation*}
z=-v
\end{equation*}%
we obtain%
\begin{equation}
2\left\vert \sigma \right\vert ^{2}\Delta v+\left\vert \nabla v\right\vert
^{2}=4\left\vert x\right\vert ^{2}+4\alpha v\text{ for }x\in \mathbb{R}^{N},
\label{las}
\end{equation}%
or equivalently%
\begin{equation}
\Delta v=\frac{4\left\vert x\right\vert ^{2}+4\alpha v-\left\vert \nabla
v\right\vert ^{2}}{2\left\vert \sigma \right\vert ^{2}}\text{ for }x\in 
\mathbb{R}^{N},  \label{las2}
\end{equation}%
Then, after changing the variable $u\left( x\right) =e^{\frac{v\left(
x\right) }{2\left\vert \sigma \right\vert ^{2}}}$, the equation (\ref{las2})
becomes%
\begin{equation*}
\left\{ 
\begin{array}{l}
\Delta u\left( x\right) =\frac{1}{\left\vert \sigma \right\vert ^{4}}%
\left\vert x\right\vert ^{2}u\left( x\right) +\frac{2\alpha }{\left\vert
\sigma \right\vert ^{2}}u\left( x\right) \ln u\left( x\right) \text{ for }%
x\in \mathbb{R}^{N}, \\ 
u\left( x\right) >0\text{ for }x\in \mathbb{R}^{N},%
\end{array}%
\right.
\end{equation*}%
which is exactly the same with the equation (\ref{11}), for%
\begin{equation*}
a\left( x\right) =\frac{1}{\left\vert \sigma \right\vert ^{4}}\left\vert
x\right\vert ^{2}\text{, }b\left( x\right) =\frac{2\alpha }{\left\vert
\sigma \right\vert ^{2}}\text{, }h\left( u\right) =u\text{, }g\left(
u\right) =u\ln u\text{ }
\end{equation*}%
and $s_{0}=1$. A recapitulation of the changes of variables and notations
are 
\begin{eqnarray*}
2\left\vert \sigma \right\vert ^{2}\ln u\left( x\right) &=&v\left( x\right)
=-z\left( x\right) , \\
H\left( s\right) &=&\ln \left( \ln s+1\right) -\ln \left( \ln u_{0}+1\right)
\Longrightarrow H^{-1}\left( s\right) =e^{e^{s+\ln \left( \ln u_{0}+1\right)
}-1}, \\
\overline{P}\left( r\right) &=&\frac{r^{4}}{4\left( N+2\right) \left\vert
\sigma \right\vert ^{4}}+\frac{\alpha r^{2}}{N\left\vert \sigma \right\vert
^{2}}, \\
\underline{P}\left( r\right) &=&\left\{ 
\begin{array}{lll}
\frac{r^{4}}{4\left\vert \sigma \right\vert ^{4}\left( N+2\right) } & \text{%
if} & u_{0}=1, \\ 
\min \left\{ 1,u_{0}\ln u_{0}\right\} \left( \frac{r^{4}}{4\left( N+2\right)
\left\vert \sigma \right\vert ^{4}}+\frac{2\alpha }{\left\vert \sigma
\right\vert ^{2}}\frac{r^{2}}{2N}\right) & \text{if} & u_{0}\neq 1.%
\end{array}%
\right.
\end{eqnarray*}%
\bigskip Then, from Theorem \ref{th1} we can see that 
\begin{equation*}
z\left( x\right) <0\text{ for }x\in \mathbb{R}^{N}\text{ and }z\left(
x\right) \rightarrow -\infty \text{ as }\left\vert x\right\vert \rightarrow
\infty .
\end{equation*}%
Let us point that, since $u\left( x_{1},...,x_{N}\right) =u\left( \left\vert
x\right\vert \right) $ is defined by 
\begin{equation*}
u\left( x\right) =u_{0}+\int_{0}^{\left\vert x\right\vert
}y^{1-N}\int_{0}^{y}t^{N-1}\left( \frac{1}{\left\vert \sigma \right\vert ^{4}%
}t^{2}u\left( t\right) +\frac{2\alpha }{\left\vert \sigma \right\vert ^{2}}%
u\left( t\right) \ln u\left( t\right) \right) dtdy,
\end{equation*}%
we have that%
\begin{equation*}
z_{x_{i}}\leq 0\text{ for all }i=1,...,N.
\end{equation*}%
for any $x=\left( x_{1},...,x_{N}\right) \in $ $\mathbb{R}_{+}^{N}$. Since
the production rate were restricted to be nonnegative, in this case 
\begin{equation*}
p_{i}^{\ast }=\max \left\{ 0,-\frac{z_{x_{i}}}{2}\right\} =-\frac{z_{x_{i}}}{%
2}\text{ for }i=1,...,N\text{.}
\end{equation*}

\begin{remark}
Using the change of variable $u\left( r\right) =e^{-\frac{z\left( r\right) }{%
2\left\vert \sigma \right\vert ^{2}}}$ we can rewrite the affirmations
A1)-A2) such:

A1)'\quad there exists $\overline{C}_{1}>u_{0}$ such that%
\begin{equation}
2\left\vert \sigma \right\vert ^{2}\ln u_{0}\leq U\left( r\right) =-z\left(
r\right) \leq 2\left\vert \sigma \right\vert ^{2}\ln \overline{C}_{1}\text{
for any }r\in \left[ 0,c_{0}\right] ,  \label{f1}
\end{equation}

A2)'\quad there exists $\overline{C}_{2}>0$ such that%
\begin{equation}
U^{\prime }\left( r\right) =-z^{\prime }\left( r\right) \leq \frac{%
2\left\vert \sigma \right\vert ^{2}}{u_{0}}\overline{C}_{2}\left( r+1\right) 
\text{ for any }r\in \left[ 0,c_{0}\right] .  \label{f2}
\end{equation}
\end{remark}

\begin{remark}
In the model problem the result of existence of solution, obtained in
Theorem \ref{th1}, holds and for the case $N\in \left\{ 1,2\right\} $.
\end{remark}

\begin{remark}
In the optimization criterion \eqref{opti} the quadratic loss function, $%
f_{1}(x)=f_{2}\left( x\right) =\left\vert x\right\vert ^{2},$ was
considered. In a future work we plan to explore optimization criterions
involving other loss functions as well.
\end{remark}

\textbf{Open problem. }Under hypotheses of the form h1) and g1) and under
some suitable conditions on $a$ and $b$ we think that there exists $l\in %
\left[ 0,s_{0}\right] $ such that the problem 
\begin{equation*}
\Delta u=a\left( x\right) h\left( u\right) +b\left( x\right) g\left(
u\right) \text{ for }x\in \mathbb{R}^{N}\text{ (}N\geq 3\text{)},
\end{equation*}%
has a unique positive solution $u\in C^{2}\left( \left[ 0,\infty \right)
\right) $ with $u\left( x\right) \rightarrow l$. Moreover, such a solution
guarantees the existence of a unique strong solution for (\ref{cons}) which
makes the candidate optimal control admissible. Moreover, we can prove some
growth estimates as in (\ref{f1}), (\ref{f2}) and then it can be proved the
transversality condition (\ref{tr}).

\textbf{Acknowledgement. }The author would like to thank Professor T. A. P.
for valuable comments and suggestions which further improved this article.

\end{document}